\documentclass[12pt,twoside,a4paper]{amsart}
\usepackage[T1]{fontenc}
\usepackage[utf8]{inputenc}
\usepackage{amsmath}
\usepackage{amsthm}
\usepackage{amsfonts}
\usepackage{bbm}
\usepackage[UKenglish]{babel}
\usepackage{enumerate}
\usepackage{bm}
\usepackage{ae}
\usepackage{amssymb}

\usepackage{hyperref}
\hypersetup{urlcolor=blue, colorlinks = true, linkcolor=black, citecolor=black}

\usepackage{easyReview}

\usepackage{mathtools}
\usepackage{wasysym}

\theoremstyle{plain}
\newtheorem{thm}{Theorem}
\newtheorem{lem}[thm]{Lemma}
\newtheorem{prop}[thm]{Proposition}

\theoremstyle{definition}

\newtheorem{example}[thm]{Example}

\theoremstyle{remark}

\numberwithin{thm}{section}

\makeatletter
\newcommand*{\house}[1]{%
	\mathord{%
		\mathpalette\@house{#1}%
	}%
}
\newcommand*{\@house}[2]{%
	\dimen@=\fontdimen8 %
	\ifx#1\scriptscriptstyle\scriptscriptfont
	\else\ifx#1\scriptstyle\scriptfont
	\else\textfont\fi\fi
	3 %
	\sbox0{%
		$#1%
		\vrule width\dimen@\relax
		\overline{%
			\kern2\dimen@
			\begingroup % to keep changes of \dimen@ in #2 local
			#2%
			\endgroup
			\kern2\dimen@
		}%
		\vrule width\dimen@\relax
		\mathsurround=1.5\dimen@ % outside margin
		$%
	}%
	\ht0=\dimexpr\ht0-\dimen@\relax
	\dp0=\dimexpr\dp0+2\dimen@\relax
	\vbox{%
		\kern\dimen@ % reinsert previously removed space
		\copy0 %
	}%
}

\newcounter{case}

\newcommand{\bC}{\mathbb{C}}
\newcommand{\bK}{\mathbb{K}}
\newcommand{\bN}{\mathbb{N}}
\newcommand{\cO}{\mathcal{O}}
\newcommand{\bQ}{\mathbb{Q}}
\newcommand{\bZ}{\mathbb{Z}}
\newcommand{\norm}{\mathcal{N}}
\newcommand{\cS}{\mathcal{S}}

\begin{document}
	\author{MATHIAS L. LAURSEN}
	
	\address{M. L. Laursen, Department of Mathematics, Aarhus University, Ny Munkegade 118,
		DK-8000 Aarhus C, Denmark}
	
	\email{mathiaslaursen@hotmail.com}
	
	\thanks{This research is supported by the Independent Research Fund Denmark}
	
	\title[Transcendence criteria for infinite products]{Transcendence criteria for infinite products of algebraic numbers}
	
	\begin{abstract}
		Using an application of Schmidt's Subspace Theorem, this paper gives new transcendence criteria for rapidly converging infinite products of algebraic numbers.
		The paper also improves existing criteria for irrationality of products and criteria for irrationality and transcendence of infinite series.
		These results generalize a classical theorem on the irrationality of infinite series due to Erd\H{o}s.
	\end{abstract}
	\maketitle
	%\textbf{\color{blue}To do: 
		%	\begin{itemize}
			%		\item \remove{Ensure that ``equation'' and ``inequality'' are used correctly.}
			%		\item \remove{Look through other comments.}
			%		\item \remove{Make sure that similar theorems are worded similarly.}
			%		\item \remove{Check spelling, word repetition etc.}
			%		\item \remove{Remove out-commented bits}
			%		\item \remove{Concider if $d_0$ should be used in Theorems \ref{thm:transc_anAlg}, \ref{thm:products general}, \ref{thm:main an algebraic} and Lemma \ref{lem:InSpan}.}
			%		\item Final: Look for unwanted line or page breaks.
			%	\end{itemize}
		%}
\section{Introduction}\label{sec:introduction}\label{sec:main results}

Proving whether a given real number is algebraic, or even rational, can be a quite frustrating endeavour.
While more than a century and a half have passed since Hermite proved that the number $e=\sum_{n=0}^{\infty} \frac{1}{n!}$ is transcendental in early 1873, it remains unsolved if the number $\sum_{n=0}^{\infty} \frac{1}{n!+1}$ is even irrational, despite what may appear as a much similar construction.
Similarly, it is well-known that the Riemann $\zeta$ function defined as $\zeta(s)=\sum_{n=1}^{\infty} n^{-s}$ for $\Re(s)>1$ is transcendental when $s$ is a positive even integer while the question of irrationality remains open when $s\ge 5$ is any fixed odd integer.
In other words, we have a multitude of interesting numbers that we know to be transcendental but where a small perturbation to the infinite series used to describe them renders even the question of irrationality exceedingly hard to settle.
Aiming away from frustrations of this kind, this paper studies irrationality and transcendence criteria that are less sensitive to such perturbations. 

Following the notions of Han\v{c}l \cite{Hancl1993,Hancl1996}, we say that a sequence $\{a_n\}_{n=1}^{\infty}$ of real or complex numbers is $\Sigma$-irrational (respectively $\Sigma$-transcenden\-tal) if the sum of the series $\sum_{n=1}^{\infty}\frac{1}{a_n c_n}$ is irrational (respectively transcendental) for any sequence $\{c_n\}_{n=1}^\infty$ of positive integers.
%Inspired by this, we say that $\{a_n\}_{n=1}^{\infty}$ is ($\Pi,K$)-irrational if the number $\prod_{n=1}^{\infty}\big(1+\frac{1}{a_nc_n}\big)$ lies outside of a given field $K$ 
Inspired by this, we say that $\{a_n\}_{n=1}^{\infty}$ is ($\Pi,K$)-irrational if the number $\prod_{n=1}^{\infty}\big(1+\frac{1}{a_nc_n}\big)$ lies outside of a given field $K$ for all sequences $\{c_n\}_{n=1}^\infty$ of positive integers. We then say that $\{a_n\}_{n=1}^{\infty}$ is $\Pi$-irrational if it is $\Pi_\bQ$-irrational and that it is $\Pi$-transcendental if it is $\Pi_{\overline\bQ}$-irrational.

The first result on $\Sigma$-irrationality was proven in 1975 by Erd\H{o}s \cite{ErdosIrrOfInfSeries}.
\begin{thm}[Erd\H{o}s]\label{thm:Erdos}
	Let $\varepsilon>0$, and let $\{a_n\}_{n=1}^\infty$ be an increasing sequence of integers such that $a_n\ge n^{1+\varepsilon}$ for all $n$.
	Suppose
	\begin{align*}
		\limsup_{n\to\infty} a_n^{2^{-n}} = \infty.
	\end{align*}
	Then the sequence $\{a_n\}_{n=1}^\infty$ is $\Sigma$-irrational.
\end{thm}
Since then, more $\Sigma$-related results have come to light, such as criteria for $\Sigma$-transcendence (starting with \cite{Hancl1996} in 1996) or $\bQ$-linear independence of the numbers $1,\sum_{n=1}^{\infty}\frac{b_{1,n}}{a_{1,n}c_n},\ldots,\sum_{n=1}^{\infty}\frac{b_{K,n}}{a_{K,n}c_n}$ (starting with \cite{Hancl1999} in 1999).

Meanwhile, the first result on $\Pi$-irrationality was not published until 2011, where Han\v{c}l and Kolouch \cite{Hancl+Kolouch2011} proved that sequences $\{a_n\}_{n=1}^\infty$ satisfying the assumptions of Theorem \ref{thm:Erdos} will also be $\Pi$-irrational.
This result was extended by Kristensen and the current author in \cite{Kristensen+Laursen2025Products} where $a_n$ are allowed to be algebraic integers from a broader family of algebraic numbers and where a lower bound on algebraic degree of the numbers $\prod_{n=1}^{\infty}(1+ a_n^{-1})$ is given.
Recall that an algebraic integer is an algebraic number whose primitive polynomial over $\bZ$ is monic.
We use $\cO_\bK$ to denote the ring of algebraic integers contained in a given field $\bK$.
When only considering the questions of $\Pi$-irrationality and $\Pi$-transcendence, this result specializes to the below theorem.
To the current author's knowledge, this is so far the only available result regarding $\Pi$-transcendence.
The notation $\house{a_n}$ denotes the maximum modulus amongst the conjugates of $a_n$.
\begin{thm}[Kristensen and Laursen]\label{thm:KL1}
	Let $\bK$ be a number field of degree $d$, let $\tilde{d},D\in\bN$ be positive integers, and consider $\varepsilon>0$ and $\alpha\in(0,1)$. 
	Let $\{a_n\}_{n=1}^\infty$ be a sequence of algebraic integers such that $n^{1+\varepsilon}< |a_n| \le |a_{n+1}|$.
	Write $\tilde\bK=\bQ(a_n:n\in\bN)$, let $\bK'$ be a field extension of degree $D$ of $\tilde\bK$, and let $\{b_n\}_{n=1}^\infty$ be a sequence of positive integers.
	Suppose $\house{a_n} b_n \le |a_n| 2^{\log_2^\alpha |a_n|}$, $\deg_{\bK} a_{n}\le \tilde d$, and $e\Re(a_n/b_n - 1/4)\ge 1/4$ for all $n$ with $e\in\{-1,1\}$ fixed and $\Re(a_n/b_n)\ne-1/2$ infinitely often.
	Then $\{a_n\}_{n=1}^\infty$ is $\Pi_{\bK'}$-irrational if $|a_n|^{D^{-n}\prod_{i=1}^{n-1}(\tilde d^i d + \tilde d)^{-1}} = \infty$ diverges in $\mathbb{R}$,
	and $\{a_n\}_{n=1}^\infty$ is $\Pi$-transcendental if for all $A>0$,
	\begin{equation*}
		\limsup_{n\to\infty} |a_n|^{A^{-n}\prod_{i=1}^{n-1}(d^i + d)^{-1}} = \infty.
	\end{equation*}
\end{thm}
In the present paper, we will improve this result in the case $\tilde d=1$, i.e., when $a_n\in \bK$ for all $n$, while providing conditions for when $\house{a_n} |b_n|$ is large and when $b_n$ is picked in $\bK$ rather than $\bQ$.
The most significant improvement lies in the transcendence criterion, where the new version allows one to stop at a finite $A$.
Results of this nature was proven for infinite series in \cite{Laursen2023subspSum} by the current author.
In that paper, the most simple criteria were found when $a_n$ is assumed to be rational while $b_n$ carries the algebraic degree of the number.
Revisiting these ideas, we will not only achieve stronger results for infinite products but also improve the theorems from \cite{Laursen2023subspSum}.
\section{Main results}
%In this paper, we make corresponding work on the $\Pi$ side.
%As it is unlikely to get $\prod_{n=1}^{N-1}(1+\frac{b_n}{a_n c_n})$ as a nice linear combination of the $x_i$, even for $a_n\in\bN$, we will not get $\Pi$-transcendence under the assumptions of Theorem \ref{thm:sums main}.
%Fortunately, Theorem \ref{thm:sums general} is more readily translated into products, though this too is slightly weakened.
As in \cite{Laursen2023subspSum}, our main results will be variations of each other.
When restricting our attention to rational numbers, we get the following irrationality and transcendence criteria.
\begin{thm}\label{thm:products over Q}
	Let $\varepsilon>0$, $0<\alpha<1$, and $\beta\in[0,\frac{\varepsilon}{1+\varepsilon})$.
	Let $\{a_n\}_{n=1}^\infty$ and $\{b_n\}_{n=1}^\infty$ be sequences of positive integers such that
	\begin{equation*}
		n^{1+\varepsilon} < a_n \le a_{n+1}
		\quad\text{and}\quad
		b_n \le a_n^\beta 2^{\log_2^\alpha a_n}.
	\end{equation*}
	Then the sequence $\{a_n/b_n\}_{n=1}^\infty$ is $\Pi$-irrational if
	\begin{equation*}
		\limsup_{n\to\infty} a_n^{(\frac{1}{1-\beta}+1)^{-n}}
		= \infty,
	\end{equation*}
	and it is $\Pi$-transcendental if
	\begin{equation*}
		\limsup_{n\to\infty} a_n^{( \frac{2+\delta}{1-\beta}+1)^{-n}}
		= \infty.
	\end{equation*}
\end{thm}

Moving on to algebraic numbers, we need to be more careful and to take into account the arithmetic properties of $a_n$ and $b_n$.
We will first consider the question of irrationality, which is also the easiest to prove.
We provide three different \emph{limsup} criteria, the latter two of which correspond to the irrationality conditions given in \cite{Laursen2023subspSum}, while the first condition is new and will be used for proving the transcendence criteria of the subsequent theorems.
While the theorem has a good number of assumptions, some of them may be skipped, depending on which \emph{limsup} condition one means to imply; thus inequality \eqref{eq:ain bounds} may be skipped when using condition \eqref{eq:irrational general}, while inequalities \eqref{eq:eta1Bound} and \eqref{eq:eta2Bound} may be skipped when using condition \eqref{eq:irrational broad} or \eqref{eq:irrational ints}.
In this theorem, and for the remainder of the current paper, $\norm:\overline{\bQ}\to\bQ$ denotes the map that sends an algebraic number to the product of its (algebraic) conjugates.
A reader familiar with field norms will notice that $\norm(a)$ is exactly the field norm associated with $\bQ(a)$, evaluated at $a$.
\begin{thm}\label{thm:irrational}
	Let $\bK$ be a number field with of degree $d\in\bN$, and consider real numbers $\alpha\in(0,1)$, $\varepsilon>0$, $\beta\in[0,\varepsilon/(1+\varepsilon))$, $y_1\ge 1$, $y_2\ge \beta$, $z_1\ge-y_2$, $z_2\ge 0$, and $e\in\{-1,1\}$.
	Let $\{a_n\}_{n=1}^\infty$ and $\{b_n\}_{n=1}^\infty$ be sequences of non-zero numbers in $\cO_{\bK}$ such that
	\begin{equation}\label{eq:|a_n| increases}
		n^{1+\varepsilon} \le |a_n| \le |a_{n+1}|,
	\end{equation}
	\begin{equation}\label{eq:bn bound}
		|b_n| < |a_n|^\beta 2^{\log_2^\alpha |a_n|},
	\end{equation}
	\begin{equation}\label{eq:ain bounds}
	\house{a_n} \le |a_n|^{y_1} 2^{\log_2^\alpha |a_n|},
	\end{equation}
	\begin{equation}\label{eq:bin bounds}
		\house{b_n} \le |a_n|^{y_2} 2^{\log_2^\alpha |a_n|},
	\end{equation}
	\begin{equation}\label{eq:eta1Bound}
		\house{a_n^{-1}} \le |a_n|^{z_1} 2^{\log_2^\alpha|a_n|},
	\end{equation}
	\begin{equation}\label{eq:eta2Bound}
		r_n \big|\norm(a_n/r_n)\big| \le |a_n|^{z_2} 2^{\log_2^\alpha |a_n|},
	\end{equation}
	and
	\begin{equation}\label{eq:real an/bn sequence}
		\Re\bigg(\frac{a_n}{b_n}\bigg) \begin{cases}
			\ge 0	&\text{if } e=1,	\\
			\le -1/2	&\text{if } e=-1,
		\end{cases}
	\end{equation}
	where each $r_n$ is a positive integer dividing $a_n$, and $\Re(a_n/b_n)\ne -1/2$ infinitely often.
	Let $d_0\in\bN$ and suppose $\deg (a_n/b_n)\ge d_0$ for all large enough $n$.
	Then the sequence $\{a_n/b_n\}_{n=1}^\infty$ is $\Pi_\bK$-irrational if 
	\begin{equation}\label{eq:irrational general}
		\limsup_{n\to\infty} |a_n|^{\left( \frac{d(y_2 + z_1 + z_2/d_0)}{1-\beta}+1 \right)^{-n}}
		= \infty,
	\end{equation}
	\begin{equation}\label{eq:irrational broad}
		\limsup_{n\to\infty} |a_n|^{\left( \frac{d(y_1+y_2)}{1-\beta}+1 \right)^{-n}}
		= \infty,
	\end{equation}
	or, in the case that $a_n\in\bZ$ or $b_n\in\bZ$ for each $n$, if
	\begin{equation}\label{eq:irrational ints}
		\limsup_{n\to\infty} |a_n|^{\left( \frac{d\max\{y_1,y_2\}}{1-\beta}+1 \right)^{-n}}
		= \infty.
	\end{equation}
\end{thm}
%\begin{rem}
%	If $c_n\le 2^{\log_2^\alpha a_n}$ for all $n$, then we would not need the additional assumption $z_2\ge 1$, thus giving us a weakened form of $\Pi_\bK$-irrationality for a slightly loosened condition.
%\end{rem}
The main novelty of this theorem lies in \emph{limsup} condition \eqref{eq:irrational general} since it grants $\Pi$-irrationality for sequences $\{a_n\}_{n=1}^\infty$ satisfying both $\deg a_n>1$ and
\begin{align*}
	\lim_{n\to\infty} |a_n|^{1/(2+\delta)^n} = 1,
\end{align*}
for all $\delta>0$, while granting a somewhat weaker result when the same is true for $\delta=0$, as seen in Example \ref{ex:irr of units}.
Meanwhile, all former $\Sigma$- or $\Pi$-irrationality statements in the literature require at least
\begin{align*}
	\lim_{n\to\infty} |a_n|^{(d+1)^{-n}} =\infty.
\end{align*}
when $\bQ(a_1,a_2,\ldots)$ is a finite field of degree $d$.

Theorem \ref{thm:irrational} also provides a stronger irrationality statement than those in \cite{Laursen2023subspSum}, where only $\Sigma_\bQ$-irrationality was proven.
This improvement, however, does not rely on any difference between products and series, and so the irrationality statements of \cite{Laursen2023subspSum} may easily be strengthened to match Theorem \ref{thm:irrational} by modifying the proofs accordingly.
An important consequence of the improved irrationality statement is that we may slack the transcendence criterion when $d>1$, giving us the below theorem.
In parallel to \cite{Laursen2023subspSum}, transcendence is proven through an application of Schmidt's Subspace Theorem, which will give us that each product $\prod_{n=1}^{\infty}(1+\frac{b_n}{c_n a_n})$ is either transcendental or contained in $\bK$, then dealing with $\bK$ through Theorem \ref{thm:irrational}. 
As with Theorem \ref{thm:irrational}, the corresponding theorem in \cite{Laursen2023subspSum} may be improved to match this result.
\begin{thm}\label{thm:products bn}
	Let $\bK$ be a number field of degree $d\in\bN$, and consider real numbers $\delta,\varepsilon>0$, $\alpha\in(0,1)$, $\beta\in[0,\frac{\varepsilon}{1+\varepsilon})$, $e\in\{-1,1\}$, and $y\ge 1$.
	Let $\{a_n\}_{n=1}^\infty$ and $\{b_n\}_{n=1}^\infty$ be a sequences of non-zero numbers from $\bZ$ and $\cO_\bK$, respectively, such that inequalities \eqref{eq:|a_n| increases}, \eqref{eq:bn bound}, \eqref{eq:bin bounds}, and \eqref{eq:real an/bn prod} are satisfied with $y_2=y$.
	Suppose that $\Re(a_n/b_n)\ne -1/2$ infinitely often and that
%	\begin{equation*}
%		|b_n| \le a_n^\beta 2^{\log_2^\alpha a_n} ,
%		\quad
%		\house{b_n} \le a_n^{y} 2^{\log_2^\alpha a_n},
%		\quad\text{and}\quad
%		e\,\Re\bigg(\frac{a_n}{b_n}+\frac{1}{2}\bigg)  > 0.
%	\end{equation*}
	\begin{equation*}
		\limsup_{n\to\infty} a_n^{\left( \frac{d y + 1 + \delta}{1-\beta}+1 \right)^{-n}}
		= \infty.
	\end{equation*}
	Then the sequence $\{a_n/b_n\}_{n=1}^\infty$ is $\Pi$-transcendental.
\end{thm}

In order to use the above mentioned application of Schmidt's Subspace Theorem, we need to write the approximants $\prod_{n=1}^{N}(1+\frac{b_n}{a_n c_n})$ as a $\bQ$-linear combination of some basis $x_1,\ldots, x_d$ of $\bK$.
If we allow $a_n\in\cO_\bK$, we then need to also consider the coordinates of $a_n^{-1}$ and the associated least common denominator, which makes both theorem and proof a bit more involved but allows us to conclude the following theorem.
Again, we get an improvement compared to \cite{Laursen2023subspSum}, and updating the proof in that paper accordingly will give a matching result for $\Sigma$-irrationality.

%Compared to Theorems \ref{thm:Erdos} and \ref{thm:KL1}, it may appear strange to see how $b_n$ carries a large part of the arithmetic information, including the algebraic degree.
%As in \cite{Laursen2023subspSum}, this comes from the proof relying on Schmidt's Subspace Theorem, to which end we need the approximants $\prod_{n=1}^{N} 1+\frac{b_n}{a_n}$ to be written as a (convenient) $\bQ$-linear combination in a fixed basis $x_1,\ldots,x_d$ of $\bK$.
%As the reader may realize, the size of the resulting coordinates have a relatively simple bound in terms of $|a_n|$ and $\house{b_n}$ when $a_n$ is rational while the same bound either increases in complexity or becomes way too rough for most applications.
%Favouring applicability over simplicity, the proof of Theorem \ref{thm:products bn} may be modified to allow irrational $a_n$ as follows.
\begin{thm}\label{thm:products general}
	Let $\bK$ be a number field with of degree $d\in\bN$, and consider real numbers $\alpha\in(0,1)$, $\delta,\varepsilon>0$, $\beta\in[0,\varepsilon/(1+\varepsilon))$, $e\in\{-1,1\}$, $y\ge\beta$, $z_1\ge-y$, and $z_2\ge 0$.
	Let $\{a_n\}_{n=1}^\infty$ and $\{b_n\}_{n=1}^\infty$ be sequences of non-zero numbers in $\cO_{\bK}$ such that inequalities \eqref{eq:|a_n| increases}, \eqref{eq:bn bound}, \eqref{eq:bin bounds}--\eqref{eq:eta2Bound}, and \eqref{eq:real an/bn prod} are satisfied with $y_2=y$, $r_n\in \bZ$, and $r_n\mid a_n$.
	Suppose that $\Re(a_n/b_n)\ne -1/2$ infinitely often and that
	\begin{equation}\label{eq:trans general}
		\limsup_{n\to\infty} |a_n|^{\left( \frac{d(y + z_1 + z_2)+ z_2 +\delta}{1-\beta}+1 \right)^{-n}}
		= \infty.
	\end{equation}
	Then the sequence $\{a_n/b_n\}_{n=1}^\infty$ is $\Pi$-transcendental.
\end{thm}
When specializing to $a_n\in\bN$, we retrieve Theorem \ref{thm:products bn}, which has the advantage of being more easily checked.
Further specializing to $b_n\in\bN$, we reach Theorem \ref{thm:products over Q}.
However, when allowing $a_n\in\cO_\bK$ and specializing to $b_n\in\bN$, we do not get the same simplification as is the case for Theorem \ref{thm:products bn}, and we will for that reason not state it as a separate theorem.

\begin{thm}\label{thm:prod only}
	Replace the assumptions $z_2\ge 1$ and \eqref{eq:real an/bn sequence} with the weaker assumptions $z_2\ge 0$ and
	\begin{equation}\label{eq:real an/bn prod}
		e\Re\bigg(\frac{a_n}{b_n} +\frac{1}{2}\bigg) \ge 0.
	\end{equation}
	Then Theorems \ref{thm:irrational}, \ref{thm:products bn}, and \ref{thm:products general} remain valid if we replace the statements of $\Pi_\bK$-irrationality and $\Pi_\bK$-transcendence with $\xi\notin\bK$ and $\xi$ being transcendental, respectively, where $\xi$ is the number $\prod_{n=1}^{\infty}(1+b_n/a_n)$.
\end{thm}

%Similarly to how condition \eqref{eq:irrational general} of Theorem Theorem \ref{thm:irrational} is unprecedentedly relaxed in terms of granting irrationality, the condition \eqref{eq:trans general} of Theorem \ref{thm:products general} can also be more relaxed than any other condition for $\Pi$- or $\Sigma$-transcendence of sequences that the current author has been able to find in the literature.
%To be precise, all these transcendence criteria require at least
%\begin{equation*}
%	\lim_{n\to\infty} |a_n|^{(3+\delta)^{-n}} = \infty,
%\end{equation*}
%while Examples \ref{ex:trans of psi} and \ref{ex:trans of units} in the below section provide sequences $\{a_n\}_{n=1}^\infty$ that satisfy the conditions of Theorem \ref{thm:products general} while 
%\begin{equation*}
%	\lim_{n\to\infty} |a_n|^{3^{-n}} = 1,
%\end{equation*}
%by having $\beta=y=z_2=0$ and $z_1<2/d$, or $\beta=y=z_1=0$ and $z_2\le1/d$, respectively.

\section{Applications}\label{sec:ex}

In the following examples, $\varphi$ is the golden ratio (the positive root of $x^2-x-1$) and $\psi$ is the supergolden ratio (the real root of $x^3-x^2-1$). Let $F_n$ and $\hat{F}_n$ be the corresponding linear recurrences, i.e., $F_1=F_2=1$, $F_{n+2}=F_n+F_{n+1}$, $\hat{F}_1= \hat{F}_2 = \hat{F}_3=1$, and $\hat{F}_{n+3}=\hat{F}_n+\hat{F}_{n+2}$.
Notice that $\varphi$ and $\psi$ are units with $\house{\varphi^{-1}}=\varphi$ and $\house{\psi^{-1}} = \psi^{1/2}$.
We will also use the notation of $\lceil a\rceil$ to denote the smallest integer $k\ge a$ when $a\in\mathbb{R}$.

If for a given example below, some of the assumptions \eqref{eq:|a_n| increases} through \eqref{eq:real an/bn prod} are not be satisfied for all of the first finitely many $n$, apply first the relevant theorem first on $\{a_{n+N}/b_{n+N}\}_{n=1}^\infty$ and then realize that the conclusion then also holds for $\{a_{n}/b_{n}\}_{n=1}^\infty$.
\begin{example}
	Let $\{h_n\}_{n=1}^\infty$ be a strictly increasing sequence of integers with $h_n\ge 3^n n$ infinitely often.
	Relying on condition \eqref{eq:irrational ints} of Theorem \ref{thm:irrational}, the sequences $\{\varphi^{h_n}\}_{n=1}^\infty$, $\{\varphi^{h_n} + b_n\}_{n=1}^\infty$, and $\{F_{h_n}/b_n\}_{n=1}^\infty$ are all $\Pi_{\bQ(\varphi)}$-irrational when $b_n\in\bQ(\varphi)$ with $0<b_n\le 2^{h_n^\alpha}$ and $\house{b_n}\le F_{h_n}$.
%	Similarly any increasing sequence combining entries from the above three sequences will also be $\Pi_{\bQ(\varphi)}$-irrational.
	The same is true if $\varphi$ is replaced by any other quadratic irrational number $x$ with $\house{x}=x>1$.
\end{example}
%Notice that while condition \eqref{eq:irrational general} would be equally applicable on either sequence $\{\varphi^{h_n}\}_{n=1}^\infty$, $\{\varphi^{h_n} + b_n\}_{n=1}^\infty$, $\{F_{h_n}/b_n\}_{n=1}^\infty$, it may not be so for combinations of $\{x^{h_n} + b_n\}_{n=1}^\infty$ and $\{\lceil x\rceil/b_n\}_{n=1}^\infty$, such as when $b_n\in\bN$ (forcing $d_0=1$) and $|\norm(x)|>x$, as is the case for, e.g., $x = \varphi + 3$ and $b_n = 1$.

The next example shows how relaxed condition \eqref{eq:irrational general} becomes when the numbers $a_n$ have the right arithmetic properties.
\begin{example}\label{ex:irr of units}
	For $i=1,2,3$, let $\{h_{i,n}\}_{n=1}^\infty$ be a strictly increasing sequence of integers with $h_{i,n}\ge (1+1/i)^n \log n$ infinitely often.
	Using condition \eqref{eq:irrational general}, Theorem \ref{thm:irrational} ensures that the sequences $\{{F}_{ h_{1,n}}\varphi^{h_{1,n}}\}_{n=1}^\infty$ and
	$\{\hat{F}_{\lceil h_{1,n}/2\rceil}\psi^{h_{1,n}}\}_{n=1}^\infty$ are $\Pi_{\bQ(\varphi)}$-irrational and $\Pi_{\bQ(\psi)}$-irrational, respectively,
	while we may use Theorem \ref{thm:prod only} to further get $\prod_{n=1}^{\infty}(1+{F}_{ h_{2,n}}^{-1}\varphi^{-h_{2,n}})\notin \bQ(\varphi)$ and $\prod_{n=1}^{\infty}1+\hat{F}_{\lceil h_{3,n}/2\rceil}^{-1}\psi^{-h_{3,n}}\notin \bQ(\psi)$.
	
	More generally, if $x>1$ is an algebraic unit with $\house{x^{-1}}=x^z$, then $\{\lceil x^{zh_{1,n}}\rceil x^{h_{1,n}}\}_{n=1}^\infty$ is $\Pi_{\bQ(x)}$-irrational while $\prod_{n=1}^{\infty}(1+\frac{1}{\lceil x^{-zh_{n}}\rceil^{-1}x^{-h_{n}}})\notin \bQ(x)$ when $\{h_n\}_{n=1}^\infty$ is a strictly increasing sequence of integers with $h_n\ge (1+Z)^{n}\log n$ infinitely often where $Z = z/(1+z)$. %Notice, however, that then $z\ge 1/(\deg x - 1)$ and, consequently, $Z\ge 1/\deg x$, so $\sup_{n\in\bN} h_n^{1/n}$ is bounded from below when $\deg x$ is fixed.
\end{example}

We then provide an example where condition \eqref{eq:irrational broad} is preferred over condition \eqref{eq:irrational general} and where we may have $a_n,b_n\notin\bZ$ infinitely often.
\begin{example}
	Let $\{h_{n}\}_{n=1}^\infty$ be a strictly increasing sequence of integers with $h_n\ge 7^n \log n$ infinitely often.
	Then condition \eqref{eq:irrational general} of Theorem \ref{thm:irrational} ensures that the sequence $\{(2^{h_n} + \sqrt[3]{2}^n )/ b_n\}_{n=1}^\infty$ is $(\Pi,\bQ(\sqrt[3]{2}))$-irrational if $b_n\in\cO_{\bQ(\sqrt[3]{2})}$ with $0<b_n<2^{h_n^\alpha}$ and $\house{b_n}\le 2^{h_n}$ for all $n$.
	The same is true when $\sqrt[3]{2}$ is replaced for any other $d$'th root of a positive integer and $h_n\ge (2d+1)^n \log n$ infinitely often. 
\end{example}

We will now give examples of $\Pi$-transcendental sequences. The first one shows how simple it is to apply Theorem \ref{thm:products bn} while the following two show how lenient equation \eqref{eq:trans general} of Theorem \ref{thm:products general} can be when the right sequences are considered.
\begin{example}
	For each $i\in\bN$, let $\{h_{i,n}\}_{n=1}^\infty$ be a strictly increasing sequence of integers with $h_{i,n}\ge (i+1/i)^n$ infinitely often.
	By Theorem \ref{thm:products bn}, $\{F_{h_{4,n}}/(1+\varphi^{-h_n})\}$ is $\Pi$-transcendental by picking $\beta=0$ and $y=1$.
	More generally, if $x>1$ is an algebraic unit of degree $d$, then the sequence $\{\lceil x^{h_{d+2, n}}\rceil /(1+x^{-h_{d+2,n}})\}_{n=1}^\infty$ is $\Pi$-transcendental.
\end{example}
\begin{example}\label{ex:trans of units}
	For $i\in\bN$, let $\{h_{i,n}\}_{n=1}^\infty$ and $\{h_{i,n}'\}_{n=1}^\infty$ be strictly increasing sequences of integers with $h_{i,n}\ge (i+1/i)^n$ and $h_{i,n}'\ge (3 - 1/i)^n$ infinitely often.
	Then Theorem \ref{thm:products general} implies that the sequences $\{{F}_{ h_{4,n}}\varphi^{h_{4,n}}\}_{n=1}^\infty$ and
	$\{\hat{F}_{\lceil h_{5,n}/2\rceil}\psi^{h_{3,n}}\}_{n=1}^\infty$ are $\Pi$-transcendental while Theorem \ref{thm:prod only} ensures that also the numbers $\prod_{n=1}^{\infty}(1+{F}_{ h_{3,n}'}^{-1}\varphi^{-h_{3,n}'})$ and $\prod_{n=1}^{\infty}1+\hat{F}_{\lceil h_{2,n}'/2\rceil}^{-1}\psi^{-h_{2,n}'}$ are transcendental.
	
	More generally, if $x>1$ is an algebraic unit with $\house{x^{-1}}=x^z$ and $\deg x = d$, then $\{\lceil x^{zh_{d+2,n}}\rceil x^{h_{d+2,n}}\}_{n=1}^\infty$ is $\Pi$-transcendental while we further have $\prod_{n=1}^{\infty}(1+\frac{1}{\lceil x^{-zh_{n}}\rceil^{-1}x^{-h_{n}}})\notin \bQ(x)$ when $\{h_n\}_{n=1}^\infty$ is a strictly increasing sequence of integers with $h_n\ge (2+dZ)^{n}$ infinitely often where $Z = z/(1+z)$.
\end{example}
\begin{example}\label{ex:trans of psi}
	Let $\{h_n\}_{n=1}^\infty$ and $\{h_n'\}_{n=1}^\infty$ be a strictly increasing sequence of integers with $h_n\ge (7-1/3)^n$ and $h_n'\ge (3-1/3)^n$ infinitely often.
	By Theorems \ref{thm:products general} and \ref{thm:prod only}, $\{\psi^{h_n}\}$ is $\Pi$-transcendental, and the number $\prod_{n=1}^{\infty}(1+\psi^{-h_n'})$ is transcendental.
	A similar argument can be made for any other real algebraic unit $x$ with $1<\house{x^{-1}}<x^{2/d}$ and assuming $h_n\ge (3-\delta)^n$ for some sufficiently small $\delta>0$.
\end{example}

\section{Preliminaries}\label{sec:prelim}
In this section, we provide useful definitions and lemmas that are either proven elsewhere or that have particularly short proofs.
	
We start by introducing two lemmas that are proven in \cite{Laursen2023subspSum}.
The first lemma uses the method of proof introduced in \cite{ErdosIrrOfInfSeries}, while the other one is an application of Schmidt's Subspace Theorem \cite{SchmidtBook}.
\begin{lem}\label{lem:ZnSmall}
	Let $\varepsilon>0$, $0<\alpha<1\le M$, and $\beta\in[0,\frac{1+\varepsilon}{\varepsilon})$.
	Let $\{a_n\}_{n=1}^\infty$ and $\{b_n\}_{n=1}^\infty$ be sequences of numbers in $\bC$ that satisfy inequalities \eqref{eq:|a_n| increases} and \eqref{eq:bn bound}.
	If 
	\begin{equation*}
		\limsup_{n\to\infty} |a_n|^{\left(\frac{M}{1-\beta}+1\right)^{-n}} = \infty,
	\end{equation*}
	then for all fixed $0<c<1$,
	\begin{equation*}
		\liminf_{N \rightarrow \infty}  \sum_{n=N+1}^\infty \left\vert \frac{b_n}{a_n}\right\vert \left(\prod_{n=1}^N |a_n|^M\right) 2^{N^2 \log_2^c |a_{N-1}|} = 0.
	\end{equation*}
\end{lem}
\begin{lem}\label{lem:appliedSubspace}
	Let $x_1,\ldots, x_d,s$ be algebraic numbers such that $s$ is $\bQ$-linearly independent of $x_1,\ldots,x_d$, and let $C,\delta>0$.
	Then the inequality
	\begin{equation}
		\label{eq:linearFormsBound}
		\Bigg|	qs - \sum_{i=1}^{d} p_i x_i	\Bigg| \prod_{i=1}^{d} \max\{1,|p_i|\} < q^{-\delta}
	\end{equation}
	has only finitely many solutions $(p_1,\ldots,p_d,q)\in\bZ^d\times\bN$ with $|p_i|\le q^C$.
\end{lem}

	The below lemma can be extracted from the proof of Lemma 15 in \cite{Kristensen+Laursen2025Products} and replaces the triangle inequality from the series setting.
	\begin{lem}\label{lem:prod->sum}
		Let $\{a_n\}_{n=1}^\infty$ be a sequence of complex numbers such that the infinite product $\prod_{n=1}^{\infty} |1+a_n|$ converges monotonously.
		Then, for all $N$,
		\begin{equation*}
			\bigg|	\prod_{n=1}^{\infty} (1+a_n)  - \prod_{n=1}^{N} (1+a_n)	\bigg| \le \max\Bigg\{1,\prod_{n=1}^{\infty} |1+a_n|\Bigg\} \sum_{n=N+1}^{\infty} |a_n|.
		\end{equation*}
	\end{lem}
		
	We now present certain notions to further describe algebraic numbers.
	Let $a$ be an algebraic number with minimal polynomial $\sum_{i=0}^{d}c_i X^i$ over the integers, with $c_d>0$.
	The leading coefficient, $c_d$, is also called the \emph{denominator} of $a$, since $c_d a$ is an algebraic integer while $c' a$ is not for any rational integer $0<c'<c_d$.
	
	By rewriting the minimal polynomial of $a$ as $c_d \prod_{i=1}^{d}(X-a_i)$ instead, 
	we define the Mahler measure as
	\begin{equation*}
		M(a) := c_d \prod_{i=1}^{d} \max\{1,|a_i|\}.
	\end{equation*}
	Surprisingly closely related to this is the Weil height, which we define as
	\begin{align*}
		H(\alpha) := \prod_{\nu\in M_\bK} \max\{1,|a|_\nu\}^{[\bK_\nu:\bQ_\nu]/[\bK:\bQ]},
	\end{align*}
	where $\bK$ is any number field containing $a$, $M_\bK$ denotes the set of places of $\bK$, $\bK_\nu$ is the local field of $\bK$ at $\nu$, and $[\bK:\bK']$ denotes the degree of a field extension $\bK\supseteq \bK'$.
	This does not depend on the choice of $\bK$ (see \cite{Waldshcmidt} for a proof).
	We will compare and estimate the house, Mahler measure, and Weil height using the following classical results.
	\begin{lem}\label{lem:heightByHouse}
		Let $a$ be an algebraic number with denominator $c_d$.
		Then
		\begin{equation*}
			H(\alpha)^d = M(\alpha) \le 
			|c_d|\max\{1,\house{\alpha}^d\}.
		\end{equation*}
	\end{lem}
	\begin{proof}
%		See \cite{Waldshcmidt}.
		The inequality is a trivial consequence of the definitions.
		For the equality, see Lemma 3.10 of \cite{Waldshcmidt}.
	\end{proof}
	\begin{lem}\label{lem:WeilHeight}
		Let $a_1,\ldots,a_n\in\overline{\bQ}$ with $a_1\ne 0$.
		Then
		\begin{gather*}
			H(a_1+\cdots+a_n)\leq nH(a_1)\cdots H(a_n),
			\quad
			H(a_1 a_2)\leq H(a_1)H(b_1),
			\\
			H\left(1/a_1\right) = H(a_1).
		\end{gather*}
	\end{lem}
	\begin{proof}
		See \cite{Waldshcmidt}.
	\end{proof}
	\begin{lem}[Liouville Inequality]
		\label{lem:LiouvilleMignotte}
		Let $\alpha$ be a non-zero algebraic number.
		Then
		\begin{align*}
			|\alpha| \geq \big(2H(\alpha)\big)^{-\deg(\alpha)}.
		\end{align*}
	\end{lem}
	\begin{proof}
		This can be extracted from \cite[Theorem A.1]{Bugeaud}.
	\end{proof}
	
	Finally, we will be using the below lemma to compare the house with a given max norm on a finite field $\bK$ seen as a $\bQ$-vector space.
	\begin{lem}\label{lem:houseOfSum}
		Let $a_1, \ldots, a_d\in \overline{\bQ}$. Then there is a constant $C_1>0$, depending only on $a_1,\ldots,a_d$, so that for any $(c_1, \ldots, c_d)\in\bQ^d$,	
		\begin{equation*}
			\house{c_1a_1 + c_2a_2 + \cdots + c_da_d} \le C_1 \max_{1\le i \le d}|c_i|.
		\end{equation*}
		If $a_1,\ldots a_d$ are linearly independent over $\bQ$, then there also is a constant $C_2>0$ such that
		\begin{equation*}
			\house{c_1a_1 + c_2a_2 + \cdots + c_da_d} \ge C_2 \max_{1\le i \le d}|c_i|.
		\end{equation*}
	\end{lem}
	\begin{proof}
		The first statement is proven in \cite{Laursen2023subspSum}, so we limit our attention to the second one.
		We make the proof by induction.
		The lemma is trivial for $d=1$, so suppose $d>1$ and that the statement holds for $d'=d-1$.
		If $c_i=0$ for any $i$, the induction is trivial, so suppose not.
		Write
		\begin{align*}
			\alpha = \frac{c_1a_1 + c_2a_2 + \cdots + c_da_d}{a_d} = c_1\frac{a_1}{a_d} + c_2 \frac{a_2}{a_d} + \cdots + c_{d-1} \frac{a_{d-1}}{a_d} + c_d.
		\end{align*}
		Since $a_1,\ldots, a_d$ are linearly independent over $\bQ$, $c_1,\ldots,c_d\ne 0$, and $d>1$, then $\alpha$ must be irrational.
		Let $\sigma_1,\ldots,\sigma_{D}$ be the distinct embeddings of $\bQ(a_1,\ldots,a_n)$ into $\overline{\bQ}$.
		
		We then have by the induction assumption that
		\begin{align*}
			\max_{1\le i\le D}\house{\alpha-\sigma_i(\alpha)} &= \max_{1\le i\le D}\house{\sum_{j=1}^{d-1} c_j \left(\frac{a_j}{a_d}- \sigma_i\bigg(\frac{a_j}{a_d}\bigg)\right)}
			\ge C_2' \max_{1\le i<d}\{|c_i|\},
		\end{align*}
		for some $C_2'>0$ that depends only on $a_1/a_d,\ldots, a_{d-1}/a_d$ and the set $\{\sigma_1,\ldots,\sigma_D\}$, which in turn only depend on $a_1,\ldots, a_{d}$.
		Pick $j$ such that $\house{\alpha-\sigma_j\alpha} = \max_{1\le i\le D}\house{\alpha-\sigma_i(\alpha)}$.
		Then
		\begin{align*}
			\max_{1\le i\le D}\house{\alpha-\sigma_i(\alpha)} = |\sigma_k(\sigma_j(\alpha)) - \sigma_k(\alpha)|,
		\end{align*}
		for a suitable $k$.
%		
%		
%		Now, the list $\alpha-g(\alpha), g(\alpha-g(\alpha)),\ldots, g(\alpha-g(\alpha))$ certainly contain all conjugates of $\alpha - g(\alpha)$, and so there is a $j\in\{0,\ldots,d-1\}$ such that $|g^j(\alpha) - g^{j+1}(\alpha)| = \house{\alpha-g(\alpha)}$.
		Noticing that $1/\house{1/a_d}$ is the modulus of the smallest conjugate of $a_d$, we then have
		\begin{align*}
			\house{\sum_{i=1}^{d}c_i a_i} &= \house{a_d\alpha}
			\ge \frac{\house{\alpha}}{\house{1/a_d}}
			\ge \frac{\max\{|\sigma_k(\sigma_j(\alpha))|, |\sigma_k(\alpha)|\}}{\house{1/a_d}}
			\\&
			\ge \frac{|\sigma_k(\sigma_j(\alpha)) - \sigma_k(\alpha)|}{2\house{1/a_d}}
			\ge \frac{C_2'}{2\house{1/a_d}} \max_{1\le i<d}\{|c_i|\},
		\end{align*}
		and the proof is complete.
	\end{proof}
	
	\section{Proof of Main Theorems}\label{sec:proof both alg}
	We will first prove the irrationality result of Theorem \ref{thm:prod only}, which we state below as a separate result for future reference.
	\begin{prop}\label{prop:irrational}
		Use the notation of Theorem \ref{thm:irrational} except that we assume $z_2\ge 0$ rather than $z_2\ge 1$.
		Suppose assumptions \eqref{eq:|a_n| increases}-\eqref{eq:eta2Bound} and \eqref{eq:real an/bn prod} are satisfied.
		Then the number $\prod_{n=1}^{\infty}(1+b_n/a_n)$ is not contained in $\bK$ if either condition \eqref{eq:irrational general} or \eqref{eq:irrational broad} is satisfied.
		If, for all $n\in\bN$, $b_n\in\bZ$ or $a_n\in\bZ$, then we may replace condition \eqref{eq:irrational broad} with condition \eqref{eq:irrational ints}.
	\end{prop}
	
	 Theorem \ref{thm:irrational}, combining ideas from \cite{Kristensen+Laursen2025Products} and \cite{Laursen2023subspSum}.
%	As it happens, the irrationality statement of Theorem \ref{thm:products general} follows from almost identical arguments, so we combine the proofs to one.
%	To keep things simple, we first state the irrationality statement of Theorem \ref{thm:products general} as a separate proposition.
%	\begin{prop}\label{prop:irr}
%		Using the notation of Theorem \ref{thm:products general}, suppose all assumptions of that theorem hold, except that the transcendence criterion \eqref{eq:trans general} does not have to be satisfied.
%		Then $\{a_n/b_n\}$ is $\Pi_\bK$-irrational.
%	\end{prop}
\begin{proof}[Proof of Theorem \ref{thm:irrational}]
	Write
	\begin{align*}
		\gamma = \prod_{n=1}^{\infty}\bigg(1+\frac{b_n}{a_n}\bigg)
		\quad\text{and}\quad
		\gamma_N = \prod_{n=1}^{N-1}\bigg(1+\frac{b_n}{a_n}\bigg),
	\end{align*}
	and assume towards contradiction that $\gamma\in \bK$.
	
	We start by boudning $H(a_n/b_n)$.
	Applying Lemmas \ref{lem:WeilHeight} and \ref{lem:heightByHouse} followed by inequalities \eqref{eq:ain bounds} and \eqref{eq:bin bounds}, we find
	\begin{align*}
		H\bigg(\frac{b_n}{a_n}\bigg)
		&= H\bigg(\frac{a_n}{b_n}\bigg)
		\le
		\begin{cases}
			\max\{|a_n|, \house{b_n}\}	&	\text{if }a_n\in\bZ, \\
			\max\{|b_n|, \house{a_n}\}	&	\text{if }b_n\in\bZ,	\\
			H(a_n)H(b_n)\le \house{a_n}\house{b_n}	&\text{regardless}
		\end{cases}
		\\&
		\le \begin{cases}
			|a_n|^{\max\{y_1,y_2\}} 2^{\log_2^\alpha a_n}	&\text{if } \forall\, m\in\bN: a_{m}\in \bZ\text{ or } b_m\in\bZ,	\\
			|a_n|^{y_1+y_2} 2^{2\log_2^\alpha a_n}	&	\text{otherwise.}
		\end{cases}
	\end{align*}
	Instead using only Lemma \ref{lem:heightByHouse}, we have
	\begin{align*}
		H\bigg(\frac{b_n}{a_n}\bigg) &\le \bigg|r_n \norm\bigg(\frac{a_n}{r_n}\bigg)\bigg|^{1/\deg (b_n/a_n)}
		\max\bigg\{1, \house{\frac{b_n}{a_n}}\bigg\}
		\\&
		\le \bigg|r_n \norm\bigg(\frac{a_n}{r_n}\bigg)\bigg|^{1/d_0}
		\max\big\{1, \house{a_n^{-1}}\house{b_n}\big\},
	\end{align*}
	so that inequalities \eqref{eq:bin bounds}, \eqref{eq:eta1Bound}, and \eqref{eq:eta2Bound} yield
	\begin{equation*}
		H\bigg(\frac{b_n}{a_n}\bigg) \le |a_n|^{y_2 + z_1 + z_2/d_0} 2^{3\log_2^\alpha |a_n|}.
	\end{equation*}
	Writing $y = \min\{y', y_2+z_1+z_2/d_0\}$ where
	\begin{align*}
		y' = \begin{cases}
			\max\{y_1,y_2\}	&\text{if for all }n\in\bN : a_n\in\bZ\text{ or } b_n\in\bZ,	\\
			y_1 + y_2	&\text{otherwise,}
		\end{cases}
	\end{align*}
	the above bounds on $H(a_n/b_n)$ combine to
	\begin{equation*}%\label{eq:H(bn/an)}
		H\bigg(\frac{b_n}{a_n}\bigg) \le |a_n|^{y} 2^{3\log_2^\alpha |a_n|}.
	\end{equation*} 
	Then, by Lemma \ref{lem:WeilHeight},
	\begin{equation}\label{eq:H(gammaN) end}
		H(\gamma-\gamma_N) \le 2H(\gamma)\prod_{n=1}^{N-1}\bigg(2 H\bigg(\frac{a_n}{b_n}\bigg)\bigg)
%		\le
%		2^N H(\gamma) \prod_{i=1}^{N-1} \big(|a_n|^y 2^{3\log_2^\alpha |a_n|}\big)
		\le 2^{2N \log_2^\alpha |a_{N-1}|} \prod_{n=1}^{N-1} |a_n|^{y},
	\end{equation}
	for large values of $N$.
	
	We now wish to apply Lemma \ref{lem:LiouvilleMignotte}.
	By inequalities \eqref{eq:|a_n| increases} and \eqref{eq:bn bound}, $|1+b_n/a_n|> 0$ for all $n$.
	Hence, inequality \eqref{eq:real an/bn prod} and the fact that $\Re(a_n/b_n)\ne-1/2$ infinitely often makes $|\gamma_N|$ monotonous with $|\gamma_N|\ne |\gamma|$, being non-decreasing for $e=1$ and non-increasing for $e=-1$.
	Hence $\gamma - \gamma_N\ne 0$ are different for all large enough $N$, and we may thus apply the lemma and then inequality \eqref{eq:H(gammaN) end} to find
	\begin{align*}
		|\gamma - \gamma_N| \ge (2H(\gamma-\gamma_N))^{-\deg(\gamma-\gamma_N)}
		\ge \Bigg(2^{3N \log_2^\alpha |a_{N-1}|} \prod_{n=1}^{N-1} |a_n|^{y}\Bigg)^{-d}.
	\end{align*}
	Hence, by Lemma \ref{lem:prod->sum},
	\begin{align*}
		\prod_{n=1}^{N-1} |a_n|^{-dy}
		&\le 
		2^{3dN \log_2^\alpha |a_{N-1}|} |\gamma - \gamma_N|
		\le 2^{4dN \log_2^\alpha |a_{N-1}|} \sum_{n=N}^{\infty}\bigg|\frac{b_n}{a_n}\bigg|,
	\end{align*}
	and thereby
	\begin{align*}
		\lim_{N \rightarrow \infty}  \sum_{n=N+1}^\infty \left\vert \frac{b_n}{a_n}\right\vert \left(\prod_{n=1}^N |a_n|^{dy}\right) 2^{N^2 \log_2^\alpha |a_{N-1}|} = \infty.
	\end{align*}
	Recalling the choice of $y$ along with the relevant assumption \eqref{eq:irrational general}, \eqref{eq:irrational broad}, or \eqref{eq:irrational ints}, this contradicts Lemma \ref{lem:ZnSmall}.
	This completes the proof.
\end{proof}
	
We are now left to prove Theorem \ref{thm:products general}.
This will follow from Theorem \ref{thm:irrational} and the below result, the proof of which takes inspiration from \cite{Laursen2023subspSum}.

\begin{thm}\label{thm:main an algebraic}
	Using the notation of Theorem \ref{thm:products general}, suppose all assumptions of that theorem, except that inequality \eqref{eq:real an/bn prod} does not have to be satisfied and the assumption $z_2\ge1$ is replaced with $z_2>0$.
	Then the number $\prod_{n=1}^{\infty}(1+b_n/a_n)$ is either transcendental or contained in  $\bK$.
\end{thm}
To prove this theorem, we will need the below lemma.
\begin{lem}\label{lem:InSpan}
	Use the notation and assumptions of Theorem \ref{thm:main an algebraic}. Write $M=d(y + z_1 + z_2) + z_2 + \delta$, and let $x_1,\ldots, x_d$ span $\cO_\bK$ as a $\bZ$-module.
		Then there exist sequences $\{p_{1,N}\}_{N=1}^\infty$, \ldots, $\{p_{d,N}\}_{N=1}^\infty$, and $\{q_N\}_{N=1}^\infty$  with $p_{i,N}\in\bZ$ and $q_N\in\bN$ with $q_N>2^N$ such that the inequalities
		\begin{equation}\label{eq:piBound general}
			|p_{i,N}| \le 2^{N \log_2^{\alpha} q_N } q_N^{1 + \frac{y+z_1}{z_2}}
%			\quad
%			\text{for all } i=1,\ldots, d
		\end{equation}
		are satisfied for each $i=1,\ldots,d$ and all sufficiently large $N$, while
		\begin{equation}\label{eq:HNS-bound general}
			\Bigg|	\sum_{n=1}^{\infty} \frac{b_n}{a_n} - \frac{\sum_{i=1}^{d} p_{i,N} x_i}{q_N}	\Bigg| < \frac{1}{2^{ dN \log_2^{\alpha} q_N } q^{M/z_2}}
		\end{equation}
		is satisfied for infinitely many $N\in\bN$.
	\end{lem}
	\begin{proof}	
		Clearly, $(x_1,\ldots,x_d)$ forms a $\bQ$-linear basis of $\bK$.
		Let $\pi_1,\ldots,\pi_d:\bK\to\bQ$ denote the associated coordinate maps.
		Notice that $\pi_i$ maps integers to integers $\pi_i(\alpha) \in\bZ$ for all $\alpha\in\cO_\bK$.

		For each $n\in\bN$, pick $\tilde{a}_n\in\bZ$ of minimal modulus such that
		\begin{align*}
			\tilde{a}_{n+1} r_{n+1}\norm(a_{n+1}/r_{n+1}) \ge \tilde{a}_n r_n\norm(a_n/r_n) \ge |a_n|^{z_2},
		\end{align*}
		noting by inequality \eqref{eq:eta2Bound} that then
		\begin{align}\label{eq:tilde a_n bound}
			\tilde{a}_n r_n\norm(a_n/r_n) \le 2|a_n|^{z_2} 2^{\log_2^\alpha |a_n|}.
		\end{align}
		Then choose
		\begin{equation*}
			q_N = \prod_{n=1}^{N-1} \Bigg(\tilde a_n r_n \norm\bigg(\frac{a_n}{r_n}\bigg)\Bigg)
			\quad \text{and}\quad
			p_{i,N} = \pi_i \Bigg(q_N \prod_{n=1}^{N-1} \bigg(1 + \frac{b_n}{a_n}\bigg)\Bigg).
		\end{equation*}
		Clearly, $q_N\in\bN$ with $q_N>2^N$ for large values of $N$, while we have $p_{i,N}\in\bZ$ since
		\begin{align*}
			\Bigg(q_N \prod_{n=1}^{N-1} \bigg(1 + \frac{b_n}{a_n}\bigg)\Bigg) = \prod_{n=1}^{N-1} \bigg(\tilde{a}_n(a_n+b_n)\frac{\norm(a_n/r_n)}{a_n/r_n}\bigg) \in\cO_{\bK}.
		\end{align*}
		Noticing		
	\begin{align*}\nonumber
		\prod_{n=1}^{N} \bigg(1+\frac{b_n}{a_n}\bigg) &
		= \sum_{\cS\subseteq \{1,\ldots,N\}}	\prod_{n\in\cS}\frac{b_n}{a_n},
	\end{align*}
		it follows from linearity of $\pi_i$ and the triangle inequality
	\begin{align*}
		\bigg|\frac{p_{i,N}}{q_N}\bigg| &=
		\Bigg|\pi_i\Bigg(\prod_{n=1}^{N-1} \bigg(1 + \frac{b_n}{a_n}\bigg)\Bigg)\Bigg|
		\le \sum_{\cS\subseteq \{1,\ldots,N\}} \Bigg|\pi_i\Bigg(\prod_{n\in\cS}\frac{b_n}{a_n}\Bigg)\Bigg|.
	\end{align*}
	By Lemma \ref{lem:houseOfSum}, there then is a $C\ge 1$ such that
	\begin{align*}
		\bigg|\frac{p_{i,N}}{q_N}\bigg| &
		\le \sum_{\cS\subseteq \{1,\ldots,N\}}	C\house{\prod_{n\in\cS}\frac{b_n}{a_n}}
		\le 2^NC
		\prod_{n=1}^{N-1}\max\big\{1, \house{b_n}\house{a_n^{-1}}\big\}.
	\end{align*}
	Hence, by inequalities \eqref{eq:bin bounds} and \eqref{eq:eta1Bound} along with the fact that $z_1\ge -y$,
	\begin{equation*}
		\bigg|\frac{p_{i,N}}{q_N}\bigg|
		\le 2^{N} C \prod_{n=1}^{N-1} |a_n|^{y + z_1} 2^{2\log_2^\alpha |a_n|}
		\le 2^{z_2^\alpha N \log_2^\alpha|a_{N-1}|} \prod_{n=1}^{N-1} |a_n|^{y + z_1}
	\end{equation*}
	for all sufficiently large $N$.
	By inequality \eqref{eq:tilde a_n bound} and the choice of $q_N$, we conclude inequality \eqref{eq:piBound general}, since
	\begin{align*}
		\bigg|\frac{p_{i,N}}{q_N}\bigg|
		&< 2^{N \log_2^\alpha(\tilde{a}_n r_n \norm(a_n/r_n))}\prod_{n=1}^{N-1} (\tilde{a}_n r_n \norm(a_n/r_n))^{\frac{y + z_1}{z_2}} %2^{\frac{2}{z_2^\alpha}\log_2^\alpha (\tilde{a}_n r_n \norm(a_n/r_n))}
		\\&\le 2^{N \log_2^\alpha q_N } q_N^{\frac{y + z_1}{z_2}}.
	\end{align*}
		
		We now just have to ensure that inequality \eqref{eq:HNS-bound general} holds for infinitely many $N$.
		The choice of $p_{i,N}$ and Lemma \ref{lem:prod->sum} show that
		\begin{align*}
			\Bigg|	\prod_{n=1}^\infty \bigg(1 + \frac{b_n}{a_n}\bigg) - \frac{\sum_{i=1}^{d} p_{i,N} x_i }{q_N}	\Bigg| &
			= \Bigg|	\prod_{n=1}^\infty \bigg(1 + \frac{b_n}{a_n}\bigg) - \prod_{n=1}^{N-1} \bigg(1 + \frac{b_n}{a_n}\bigg) 	\Bigg|
			\\&
			\le C'\sum_{n=N}^\infty \bigg|\frac{b_n}{a_n}\bigg|
		\end{align*}
		holds for all sufficiently large $N$ and a suitable $C'\ge 1$.
		By inequality \eqref{eq:trans general} and the choice of $M$, Lemma \ref{lem:ZnSmall} then implies that 
		\begin{align*}
			\Bigg|	\prod_{n=1}^\infty \bigg(1 + \frac{b_n}{a_n}\bigg) - \frac{\sum_{i=1}^{d} p_{i,N} x_i }{q_N}	\Bigg| &< \left(\prod_{n=1}^N |a_n|^{-M}\right) 2^{-N^2 \log_2^\alpha |a_{N-1}|}
		\end{align*}
		for infinitely many $N$.
		For these $N$, the choices of $\tilde{a}_n$ and $q_N$ now allow us to conclude inequation \eqref{eq:HNS-bound general} by calculating
		\begin{align*}
			\Bigg|	\prod_{n=1}^\infty \bigg(1 + \frac{b_n}{a_n}\bigg) - \frac{\sum_{i=1}^{d} p_{i,N} x_i }{q_N}	\Bigg| &< \frac{q_N^{-M/z_2} }{2^{2d N^{2-\alpha} \log_2^\alpha ((r_{N-1}\tilde{a}_{N-1}\norm(a_n/r_n))^N)/z_1^\alpha}} 
			\\&
			\le \frac{1}{2^{dN \log_2^\alpha q_N} q_N^{M/z_2} }.
		\end{align*}
		This completes the proof.
	\end{proof}
	
	\begin{proof}[Proof of Theorem \ref{thm:main an algebraic}]
		By Lemma \ref{lem:InSpan}, infinitely many $(p_1,\ldots,p_d,q)\in\bZ^d\times\bN$ satisfy both inequalities \eqref{eq:piBound general} and \eqref{eq:HNS-bound general} with $|p_i|<q^C$ for some fixed $C>0$ that does not depend on $(p_1,\ldots,p_d,q)$.
		We first rewrite inequality \eqref{eq:HNS-bound general} by recalling $M=d(y + z_1 + z_2) + z_2 + \delta$,
		\begin{align*}
			\Bigg|	q\prod_{n=1}^{\infty}  \bigg(1+\frac{b_n}{a_n}\bigg) - \sum_{i=1}^{d} p_i x_i	\Bigg|\prod_{i=1}^{d} \Big(q^{1+\frac{d(y + z_1 + z_2)}{z_2}}2^{N \log_2^{\alpha} q } \Big)< q^{-\delta/z_2}.
		\end{align*}
		It then follows from the inequalities of \eqref{eq:piBound general} that
		\begin{align*}
			\Bigg|	q \prod_{n=1}^{\infty} \bigg(1+\frac{b_n}{a_n}\bigg) - \sum_{i=1}^{d} p_i x_i	\Bigg|\prod_{i=1}^{d} \max\{1,|p_i|\} 
			< q^{-\delta/z_2}.
		\end{align*}
		Lemma \ref{lem:appliedSubspace} now implies that $\prod_{n=1}^{\infty} \big(1+b_n/a_n\big)$ cannot both be algebraic and $\bQ$-linearly independent of $x_1,\ldots,x_d$.
		In other words, the number $\prod_{n=1}^{\infty} \big(1+b_n/a_n\big)$ is either contained in $\bK$ or transcendental, and the proof is complete.
	\end{proof}
\begin{proof}[Proof of Theorem \ref{thm:prod only}]
	Replace $\delta$ and $z_2$ with $\delta'=\delta/(d+2)$ and $z_2'=z_2 + \delta'$, respectively. Then the statement follows from Proposition \ref{prop:irrational} and Theorem \ref{thm:main an algebraic}.
\end{proof}
\begin{proof}[Proof of Theorems \ref{thm:irrational} and \ref{thm:products general}]
	Let $\{c_n\}_{n=1}^\infty$ be a sequence of positive integers.
	Replacing $a_n$ and $r_n$ by $a_nc_n$ and $c_nr_n$, respectively, then invalidates neither $a_n\ge n^{1+\varepsilon}$ nor any of the inequalities \eqref{eq:bn bound}-\eqref{eq:real an/bn sequence} since we clearly have $\varepsilon,\beta,y_2\ge 0$, $y_1,z_2\ge 1$, and $z_1\ge -1$.
	Rearranging the terms of $\{a_nc_n/b_n\}_{n=1}^\infty$ so that $|a_nc_{n}|$ becomes non-decreasing then allows us to apply Theorem \ref{thm:prod only} on the number $\prod_{n=1}^{\infty} \big(1+\frac{b_n}{a_n c_n}\big)$.
	This completes the proof.
\end{proof}

\textbf{Acknowledgements.} 
I thank the Independent Research Fund Denmark (Grant ref. 1026-00081B) for funding the research, and I thank my supervisor Simon Kristensen for proofreading the paper.
Finally, I thank the reviewer for providing helpful comments.

\providecommand{\bysame}{\leavevmode\hbox to3em{\hrulefill}\thinspace}
\providecommand{\MR}{\relax\ifhmode\unskip\space\fi MR }
\providecommand{\MRhref}[2]{%
	\href{http://www.ams.org/mathscinet-getitem?mr=#1}{#2}
}
\providecommand{\href}[2]{#2}

\end{document}